\input amstex
\documentstyle{amsppt}
\input bull-ppt
\keyedby{bull397e/lbd}
\NoBlackBoxes
\input BoxedEPS.tex
\input BoxedEPS.cfg
\HideDisplacementBoxes
\define\e{\epsilon} \define\s{\sigma} \define\Bd{\partial}
     \define\C{\Bbb C} 
\redefine\D{\Delta} 
\redefine\P{\Bbb P} 
\define\Pscr{\scr P}
\define\R{\Bbb R} 
\define\Sscr{\scr S}
     \define\Z{\Bbb Z} \define\rank{\operatorname{rank}}
     \hyphenation{Worces-ter Al-ex-and-er Kron-heim-er}
\def\twobars#1#2#3#4{\vcenter{\hrule height.#1pt width#2pt 
                               \vskip#3pt 
                               \hrule height.#1pt width#2pt
                               \vskip#4pt}}
     \def\stroke#1#2#3{\vrule height#1pt width.#2pt 
depth#3pt}
     
\def\connsum{\,\twobars4{1.3}31\stroke831\twobars4{1.8}31
                  \stroke831\twobars4{1.3}31\thinspace}
\topmatter
\cvol{29}
\cvolyear{1993}
\cmonth{July}
\cyear{1993}
\cvolno{1}
\cpgs{51-59}
\ratitle
\title 
Quasipositivity as an obstruction to sliceness 
\endtitle
\author 
Lee Rudolph 
\endauthor
\address 
Department of Mathematics and Computer Science, Clark 
University,
Worcester, Massachusetts 01610 
\endaddress
\ml
lrudolph\@clarku.bitnet 
\endml
\date June 10, 1992 and, in revised form, November 1, 
1992\enddate
\abstract
For an oriented link $L \subset S^3 = \Bd\!D^4$, let 
$\chi_s(L)$ be the 
greatest Euler characteristic $\chi(F)$ of an oriented 
2-manifold $F$ 
(without closed components) smoothly embedded in $D^4$ 
with boundary $L$.  
A knot $K$ is {\it slice} if $\chi_s(K)=1$.  Realize $D^4$ 
in $\C^2$ as 
$\{(z,w):|z|^2+|w|^2\le1\}$.  It has been conjectured 
that, if $V$ is a 
nonsingular complex plane curve transverse to $S^3$, then 
$\chi_s(V\cap S^3)=\chi(V\cap D^4)$.  Kronheimer and 
Mrowka have
proved this conjecture in the case that $V\cap D^4$ is the 
Milnor fiber 
of a singularity.  I explain how this seemingly special 
case implies both
the general case and the ``slice-Bennequin inequality'' 
for braids.  
As applications, I show that various knots are not slice 
(e.g., pretzel 
knots like $\Pscr(-3,5,7)$; all knots obtained from a 
positive trefoil 
$O\{2,3\}$ by iterated untwisted positive doubling).  As a 
sidelight, I 
give an optimal counterexample to the ``topologically 
locally-flat Thom 
conjecture''.  
\endabstract
\endtopmatter

\document
\heading 1. A brief history of sliceness
\endheading
A {\it link} is a compact 1-manifold without boundary $L$ 
(i\.e\., 
finite union of simple closed curves) smoothly embedded in 
the 3-sphere
$S^3$; a {\it knot} is a link with one component.  If 
$S^3$ is realized in
$\R^4$ as, say, the unit sphere, then a natural way to 
construct links 
is to intersect suitable two-dimensional subsets 
$X\subset\R^4$ with $S^3$;
one may then ask how constraints on $X$ are reflected in 
constraints on the
link $X\cap S^3$.

For instance, Fox and Milnor (c\. 1960) considered, in 
effect, the 
case that $X$ is a smooth 2-sphere intersecting $S^3$ 
transversally; at 
Moise's suggestion, Fox \cite{5} adopted the adjective 
{\it slice} to 
describe the knots and links $X\cap S^3$ so constructed.  
Fox and Milnor 
\cite{6} gave a criterion for a knot $K$ to be slice: its 
{\it Alexander 
polynomial} $\D_K(t)\in\Z[t,t^{-1}]$ must have the form 
$F(t)F(t^{-1})$.  This shows that, for instance, the two 
trefoil knots 
$O\{2,\pm3\}$ are not slice (since 
$\D_{O\{2,\pm3\}}=t^{-1}-1+t$ is not 
of the form $F(t)F(t^{-1})$), but it says nothing about 
the two granny knots 
$O\{2,3\}\connsum O\{2,3\}$, $O\{2,-3\}\connsum O\{2,-3\}$ 
(indeed, both 
granny knots share the Alexander polynomial $(t^{-1}-1+
t)^2$ with the square 
knot $O\{2,3\}\connsum O\{2,-3\}$, which {\it is} slice), 
and Fox could only 
aver \cite{5} that ``it is highly improbable that the 
granny knot is a 
slice knot.''
 
By the end of the 1960s, several mathematicians
[30, 14, 31] had found invariants which could be applied 
to show that, 
for instance, the granny knots are not slice.  For any 
knot $K$, all these 
invariants (signatures of various families of hermitian 
forms), as well as 
$\D_K(t)$, can be calculated from the {\it Seifert 
pairing} $\theta_F:H_1(F,\Z)\times H_1(F,\Z)\to \Z$ 
determined 
by any {\it Seifert surface} $F$ for $K$ (i\.e\., a 
smooth, oriented, 
2-submanifold-with-boundary $F\subset S^3$ without closed 
components, with 
$K=\Bd F$).  In particular, if $K$ is slice, then (for any 
$F$) there is a 
subgroup $N\subset H_1(F,\Z)$ with 
$\rank(N)=\frac{1}{2}\rank(H_1(F,\Z))$ 
on which $\theta_F$ vanishes identically.  Any knot for 
which such a subgroup 
exists is called {\it algebraically slice} (briefly, 
$A$-slice).  
Levine showed \cite{13} that in higher odd dimensions, 
$A$-slice knots are 
slice.  Whether this were true for 
knots in $S^3$ was unknown until 1975, when Casson and 
Gordon [1, 2] developed ``second-order'' obstructions to 
sliceness (again using signatures, but of more subtly 
constructed forms that 
are not determined just by a Seifert pairing) and used 
them to show that 
many $A$-slice knots are not slice.  Their methods were 
powerless, however, to 
prove nonsliceness of any knot $K$ with $\D_K(t)=1$ (such 
a knot 
is $A$-slice, as observed by L\. Taylor, cf\. 
\cite{9\,(1978), problem 1.36}).

The subject took surprising turns in the 1980s after 
Donaldson and 
Freedman revolutionized the theory of 4-manifolds and, not 
so incidentally,
the theory of knots and links in $S^3$.  In fact, let 
$X\subset\R^4$ now be 
a 2-sphere, still transverse to $S^3$, which is, however, 
assumed no longer 
smooth but merely {\it topologically locally-flat} (i.e., 
in local 
$\scr C^0$ charts it looks like $\R^2\subset\R^4$); then 
the link 
$X\cap S^3$ is {\it topologically locally-flatly slice} 
(briefly, $T$-slice).  
$T$-slice implies $A$-slice. Freedman \cite{8} proved that 
any knot $K$ 
with $\D_K(t)=1$ (e.g., any {\it untwisted double}) is 
$T$-slice.  
Nonsliceness results flowed from Donaldson's restrictions 
on 
intersection forms of smooth, as distinct from topological, 
4-manifolds: Casson proved the existence of a nonslice 
knot $K$ 
with $\D_K(t)=1$ (cf\. \cite{9\,(1984), problem 1.36}); 
Akbulut gave an explicit example of such a knot, 
the untwisted positive double $D(O\{2,3\},0,+)$ (cf\. 
\cite{3});
Cochran and Gompf \cite{3} found large classes of knots 
$K$ such
that $D(K,0,+)$ is not slice; and Yu \cite{32}, building 
on work of Fintushel 
and Stern, found many $A$-slice Montesinos knots which are 
not slice.

In \S4 I give many examples of nonslice knots: for 
example---recovering
some of Fintushel and Stern's results---all 
pretzel knots $\Pscr(p,q,r) \ne O$ with Alexander 
polynomial 1,
and---considerably generalizing 
\cite{3,~Corollary~3.2}---all iterated untwisted positive 
doubles of any 
knot $K \ne O$ which is a closed positive braid.  
The method in each case is to show that the knot under 
consideration is 
{\it strongly quasipositive}, then to use the fact that
a strongly quasipositive knot $K\ne O$ is not slice, which 
follows from
a corollary to a recent result of Kronheimer and Mrowka 
\cite{10}.  In \S3
I state their result and establish that corollary, as well 
as a superficially
stronger (actually equivalent) corollary, the 
``slice-Bennequin inequality''
for braids.  Section 2 is preliminary material on 
quasipositivity, etc. 
Section
5 is a 
sidelight, using an example from \S2 to produce a 
topologically locally-flat 
surface in $\C\P^2$, of algebraic and geometric degree 5, 
with genus 
$5=\frac{1}{2}(5-1)\times (5-2)-1$: this is an optimal 
counterexample to the 
``topologically locally-flat Thom conjecture''.  

\rem{Remarks} (1) Note that it is not Kronheimer and 
Mrowka's machinery, but
``only'' their (spectacular) result which is used.  In 
particular, one can
understand the present note while staying totally 
disengaged from gauge theory.
\par
(2) Although Kronheimer and Mrowka in \cite{10} do not 
discuss the
slice-
Bennequin inequality, they do draw explicit attention to a 
(strictly
weaker) corollary of their main result, namely, the 
affirmative answer to the
``question of Milnor'' \cite{15} on the unknotting number 
of a link of a
singularity.  The nonsliceness results of the present 
paper have nothing to do
with unknotting number.
\par
(3) W\. M\. Menasco has recently announced a proof of
the unknotting result which, in marked contrast to that in 
\cite{10}, uses
purely three-dimensional techniques (somewhat in the style
of \cite{0}); should such techniques someday be used 
successfully to establish 
the slice-Bennequin inequality, then the present 
nonsliceness results will
have a purely three-dimensional proof as well.
 \endrem
\heading 2. Quasipositivity
\endheading
\subheading
{Transverse $\C$-links and quasipositive Seifert surfaces}
When constructing links as intersections $X\cap S^3$, 
instead of restricting
the topological type of $X$ as in \S1, one might restrict 
the nature of
the embedding $X\hookrightarrow\R^4$.  In particular, if 
$\R^4$ is identified 
with $\C^2\supset S^3:=\{(z,w):|z|^2+|w|^2=1\}$ and $X$ is 
required to be 
a complex plane curve, then one can obtain many 
interesting links.

\dfn{Definitions} A {\it complex plane curve} is any set 
$V_f:=f^{-1}(0)\subset\C^2$, where $f(z,w)\in\C[z,w]$ is 
nonconstant; 
$V_f$ is a smooth, oriented 2-submanifold of $\C^2$ except 
at a finite 
set $\Sscr(V_f)\subset V_f$ of singularities.  If $V_f$ is 
transverse to 
$S^3$, then the oriented link $K_f:=V_f\cap S^3$ is a 
{\it transverse $\C$-link} [22, 29].
\enddfn

\dfn{Examples} Replacing $S^3$ by a round sphere of 
sufficiently small 
radius centered at a point of $\Sscr(V_f)$, one sees that 
any 
{\it link of a singularity} of a complex plane curve is a 
transverse $\C$-link;
replacing $S^3$ by a round sphere of sufficiently large 
radius, one sees 
that any {\it link at infinity} of a complex plane curve 
is a transverse 
$\C$-link.  
\enddfn

Links of singularities and links at infinity, though very 
interesting
(cf\. [15, 11, 4, 23, 17], etc.), are
highly atypical transverse $\C$-links (for instance, while 
the unknot $O$ is 
the only slice knot which is a link of a singularity 
\cite{11} or a link 
at infinity \cite{23}, many nontrivial slice knots are 
transverse 
$\C$-links \cite{19}).  A much broader class of transverse 
$\C$-links is 
easily defined using braid theory.

\dfn{Definitions} In the {\it $n$-string braid group} 
$$
B_n:=\text{gp}\left(\s_i, 1\leq i\leq n-1\left|
\matrix\format\l\\
[\s_i, \s_j]=\s_j^{-1}\s_i,\\
[\s_i, \s_j]=1,\endmatrix\quad
{{|i-j|=1}\atop{|i-j|\not=1}}\right.  \right), 
$$
a {\it positive band} is any conjugate $w\s_i w^{-1}$ 
($w\in B_n,  
1\leq i\leq n-1)$; a {\it positive embedded band} is one 
of the positive bands 
$\s_{i,j}:=(\s_i\dotsm\s_{j-2})\s_{j-1}(\s_i\dotsm%
\s_{j-2})^{-1}$, 
$1\leq i<j\leq n$ (e.g., each {\it standard generator} 
$\s_i=\s_{i,i+1}$
is a positive embedded band).  A {\it {\rm(}strongly{\rm)} 
quasipositive braid}
is any product of positive (embedded) bands (e.g., a {\it 
positive braid},
that is, a product of standard generators, is strongly 
quasipositive).  
A {\it {\rm(}strongly{\rm)} 
quasipositive} oriented link is one which can be realized 
as the closure of a 
(strongly) quasipositive braid.  Up to ambient isotopy, 
every quasipositive 
link is a transverse $\C$-link \cite{19}.
\enddfn
\dfn{Question}  Is every transverse $\C$-link quasipositive?
\enddfn
\rem{Remarks} (1)~There are non-quasipositive knots, for 
example,
the figure-8.
This follows, for instance, from a result of Morton 
\cite{16} 
and Franks and Williams \cite{7} about the {\it oriented 
link polynomial} of 
a closed braid (cf\. \cite{26}).  (Note, however, that 
every Alexander 
polynomial, and indeed every Seifert pairing, can be 
realized by a 
quasipositive knot or link \cite{21}.) 
\par
(2)~There are knots which are not 
transverse $\C$-links; the figure-8 is again an example.  
Biding an affirmative answer to the above question, I know 
of no way to show 
this without using the methods of the present paper.
\endrem

Any specific expression of a quasipositive braid as a 
product of positive 
bands, $\beta=
w_1^{\phantom1}\s_{i_1}w_1^{-1} 
w_2^{\phantom1}\s_{i_2}w_2^{-1}\dotsm%
w_k^{\phantom1}\s_{i_k}w_k^{-1}\in B_n$,
gives a recipe for constructing a {\it quasipositive 
braided Seifert ribbon}
$S(w_1^{\phantom1}\s_{i_1}w_1^{-1},\dots 
,w_k^{\phantom1}\s_{i_k}w_k^{-1})
\subset D^4$, that is, a
smooth surface (actually ``ribbon-embedded'', a refinement 
we can ignore)
bounded by the closed braid $\widehat\beta$.  
The isotopy carrying $\widehat\beta$ onto a transverse 
$\C$-link $K_f$ 
can be chosen to carry 
$S(w_1^{\phantom1}\s_{i_1}w_1^{-1},\dots,%
w_k^{\phantom1}\s_{i_k}w_k^{-1})$ onto
the (nonsingular) piece of complex plane curve $V_f\cap 
D^4$.  The Euler
characteristic of 
$S(w^{\phantom1}_1\s_{i_1}w_1^{-1},\dots 
,w_k^{\phantom1}\s_{i_k}w_k^{-1})$ 
is $n-k$.  If 
$\beta=\s_{i_1,j_1}\s_{i_2,j_2}\dotsm\s_{i_k,j_k}$ is 
strongly 
quasipositive, then 
$S(\s_{i_1,j_1},\s_{i_2,j_2},\dots,\s_{i_k,j_k})\subset 
D^4$ is the
``push-in'' of a {\it quasipositive braided Seifert 
surface}, abusively
indicated by the same notation; Figs\.~1 and~2 give a 
sufficient idea of the 
construction.   

A Seifert surface is {\it quasipositive} if it is ambient 
isotopic to a 
quasipositive braided Seifert surface.  (See \cite{20--22} 
for more on 
braided surfaces and quasipositivity.)  A subset of a 
surface is {\it full} 
if no component of its complement is contractible.  

\proclaim{Theorem {\rm\cite{18}}} A full subsurface of a 
quasipositive Seifert 
surface is quasipositive.
\endproclaim

\subheading
{Plumbing; quasipositive doubles }
For $K$ a knot, $\tau\in\Z$, let $A(K,\tau) \subset S^3$ 
be an 
{\it annulus of type $K$ with $\tau$ twists}; that is, 
$K\subset\Bd A(K,\tau)$
and $\theta_{A(K,\tau)}$ has matrix $\left( \tau \right)$.  
Let $A(K,\tau)*A(O,\pm 1)$ be a Seifert surface formed 
by {\it plumbing} $A(O,\pm 1)$ to $A(K,\tau)$; that is, 
there is a 3-cell $B\subset S^3$ such that 
$A(K,\tau)\subset B$, 
$A(O,\pm 1)\subset S^3\setminus\operatorname{Int}B$, and
$A(K,\tau)\cap A(O,\pm 1)=A(K,\tau)\cap\Bd B=A(O,\pm 
1)\cap\Bd B$ is
a quadrilateral 2-cell whose sides are, in order, 
contained in 
alternate components of $\Bd A(K,\tau)$ and $\Bd A(O,\pm 
1)$.
The knot $D(K,\tau,\pm):=\Bd (A(K,\tau)*A(O,\mp 1))$ is the 
{\it $\tau$-twisted positive} (resp. {\it negative}) 
double of $K$.  
A matrix for $\theta_{D(K,\tau,\pm)}$ is 
$\binom{\tau\;\phantom{\mp}\!1}{0\;\mp 1}$,
so $\D_{D(K,\tau,\pm)}(t)= 1\mp\tau(t-2+t^{-1})$, and 
$D(K,0,\pm)$ is $A$-slice for any $K$.
\proclaim
{Lemma 1} If $K\ne O$ is strongly quasipositive, then 
$A(K,0)$ is 
quasipositive.
\endproclaim

\midinsert
\vskip.625in
\centerline{\BoxedEPSF{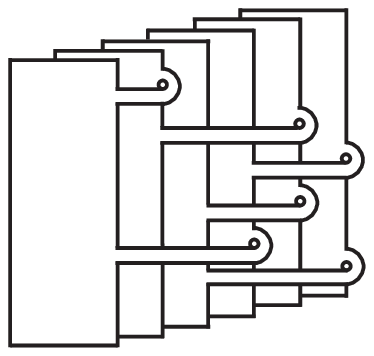 scaled 827}}
\bigskip
\centerline{{\smc Fig.\ 1}
$S(\s_{3,6},\s_{1,4},\s_{3,5},\s_{4,6},\s_{2,5},\s_1)$.}
\vskip2pc
\centerline{\BoxedEPSF{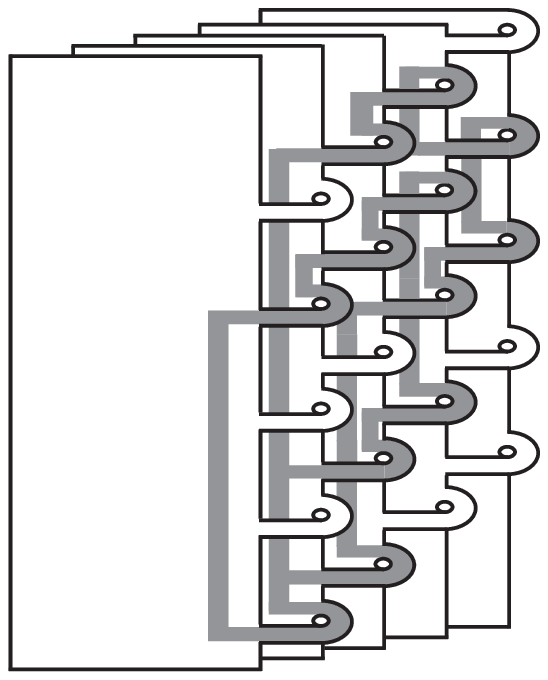 scaled 425}}
\bigskip
\centerline{{\smc Fig.\ 2}
$F(-3,5,7)$ on the Seifert surface of $O\{5,5\}$.}
\endinsert

\demo
{Proof}  This follows from the last theorem; for a collar 
of the boundary of 
a quasipositive Seifert surface $F\ne D^2$ bounded by $K$ 
is an annulus 
$A(K,0)$, and full.\qed
\enddemo
\dfn{Example} $O\{2,3\}=\Bd S(\s_1,\s_1,\s_1)$; 
$A(O\{2,3\},0)$ is 
isotopic to the quasipositive braided surface 
$S(\s_{3,6},\s_{1,4},\s_{3,5},\s_{4,6},\s_{2,5},\s_1)$ 
pictured in
Fig\.~1.
\enddfn
\proclaim
{Lemma 2} 
If the knot $K\ne O$ is strongly quasipositive, then 
$D(K,0,+)$ is strongly
quasipositive, being the boundary
of a quasipositive braided Seifert surface of Euler 
characteristic $-1$.
\endproclaim
\demo
{Proof} This follows from Lemma 1 and a theorem in 
\cite{25}: for any Seifert 
surface $S$, annulus $A$, and proper arc $\alpha \subset 
S$, the plumbed 
surface $S*_{\alpha}A$ is quasipositive if 
both $S$ and $A$ are quasipositive.  A proof in the 
present case, where 
$S$ is itself an annulus, $\alpha$ is a transverse arc of 
$S$, and 
$A=A(O,-1)$, was given in \cite{22}; the reader can readily 
recreate it after comparing the following example to the 
preceding one.\qed
\enddemo
\dfn{Example} $D(O\{2,3\},0,+)=
\Bd S(\s_6,\s_{3,6},\s_6,\s_{1,4},\s_{3,5},\s_{4,6},
\s_{2,5},\s_1)$.
\enddfn

\subheading
{Quasipositive pretzels }
Let $p, q, r \in \Z$.  A diagram for the {\it pretzel} link 
$\Pscr(p,q,r)$ is obtained from a braid diagram for 
$\beta_{p,q,r}:=\s_1^{-p}\s_3^{-q}\s_5^{-r} \in B_6$ 
by forming the {\it plat} of $\beta_{p,q,r}$ (using the 
pairing $(16)(23)(45)$
at top and bottom).  If $p, q, r$ are all 
odd, then $\Pscr(p,q,r)$ is a knot, and (once it is 
oriented) the obvious 
surface $F(p,q,r)$ that it bounds (two $0$-handles 
attached by three 
$1$-handles) is a Seifert surface.

\dfn{Example}
$\Pscr(1,1,1)=O\{2,3\}$; $F(1,1,1)=S(\s_1,\s_1,\s_1)$ (up to
ambient isotopy).
\enddfn

\proclaim{Lemma 3} For $p, q, r$ all odd{\rm,} $F(p,q,r)$ 
is quasipositive 
iff 
$$
\min\{p+q,p+r,q+r\}>0.
\tag *
$$
\endproclaim
\demo{Proof} 
For
$-\tau\in\{p+q,p+r,q+r\}$, $F(p,q,r)$ contains $A(0,\tau)$ 
as a full subsurface
(omit each $1$-handle in turn).
It is proved in \cite{29} that $A(O,\tau)$ is 
quasipositive iff $\tau<0$; 
therefore,
by the theorem of \cite{18} quoted above, if $F(p,q,r)$ is 
quasipositive,
then $(*)$ is true.  Conversely, if $(*)$ is true, then 
either $\min\{p,q,r\}>0$, or exactly one of $p,q,r$ is 
negative and it is of
strictly smaller absolute value than the other two.  In 
the first case,
$F(p,q,r)$ is obtained (up to ambient isotopy) from the 
quasipositive 
Seifert surface
$S(\s_1,\s_1,\s_1)$ by applying nonpositive twists to
the three $1$-handles, so, according to \cite{21} (or 
\cite{22}), 
$F(p,q,r)$ is quasipositive; a similar, only slightly less 
straightforward, 
twisting argument applies in the second case.\qed
\enddemo
\dfn{Example} $F(-3,5,7)$ is ambient isotopic to 
$$S(\s_1,\s_2,\s_{2,4},\s_{3,6},\s_{1,4},\s_5,\s_{2,5}).$$
\enddfn

\heading 3. Kronheimer-Mrowka Theorem; ``slice-Bennequin 
inequality''
\endheading
If $L\subset S^3$ is an oriented link, let $\chi_s(L)$ be 
the greatest Euler 
characteristic $\chi(F)$ of an oriented 2-manifold $F$ 
(without closed 
components) smoothly embedded in $D^4$ with boundary $L$; 
so, for a knot $K$,
$\chi_s(K)=1$ iff $K$ is slice.

If $K_f\subset S_\e^3$ is the link of the singularity 
$(0,0)\in\Sscr(V_f)$, then its {\it Milnor fiber} \cite{15}
is the nonsingular piece of complex plane curve 
$V_{f-\delta}\cap D_\e^4$ (for any sufficiently small 
$\delta>0$);
of course, $K_f$ is isotopic to $K_{f-\delta}=\Bd 
V_{f-\delta}\cap D_\e^4$.
The following is a restatement of 
\cite{10,~Corollary~1\.3} in the
present terminology.  

\proclaim{Kronheimer-Mrowka Theorem} If $K_f$ is the link 
of a singularity{\rm,} 
then $\chi_s(K_f)$ is the Euler characteristic of its 
Milnor fiber.
\endproclaim

This is a special case of the next proposition, which, 
however, it implies!

\proclaim{Proposition} If $K_f\subset S^3$ is a transverse 
$\C$-link and 
$\Sscr(V_f)\cap D^4=\emptyset${\rm,} then 
$\chi_s(K_f)\allowmathbreak=\allowmathbreak\chi(V_f\cap 
D^4)$.
\endproclaim

\demo{Proof} Without loss of generality (after perturbing 
$f$ slightly) 
we may assume that the projective completion 
$\Gamma\subset\C\P^2$ of $V_f$ 
in $\C\P^2\supset\C^2$ is nonsingular and transverse to 
the line at infinity.
Then the link at infinity of $V_f$ is isotopic to 
$O\{d,d\}$, $d=\deg\Gamma$.  
Assuming $\chi_s(K_f) > \chi(V_f\cap D^4)$, we would then 
also have 
$\chi_s(O\{d,d\}) > \chi(V_f)$.  Yet $O\{d,d\}$ is also 
a link of a singularity (namely, $z^d+w^d$ at the origin), 
and the
interior of its Milnor fiber is diffeomorphic to $V_f$, so 
our assumption
is inconsistent with the Kronheimer-Mrowka Theorem. \qed
\enddemo

\proclaim{Corollary} If 
$\beta=w_1^{\phantom1}\s_{i_1}w_1^{-1}\dotsm%
w_k^{\phantom1}\s_{i_k}w_k^{-1}\in B_n$ is 
quasipositive{\rm,}
then $\chi_s(\widehat\beta)=n-k$. 
\qed
\endproclaim

This corollary---in fact, its special case that 
a strongly quasipositive knot $K\ne O$ is not 
slice---suffices 
for the applications in \S4.  It is easy, however, to 
go further.  Let $e:B_n\to\Bbb Z$ be abelianization 
(exponent 
sum with respect to the standard generators $\s_i$).

\proclaim{Slice-Bennequin Inequality} For every $n,$ for 
every $\beta\in B_n,$ 
$\chi_s(\widehat\beta) \le n-e(\beta)$.
\endproclaim

\demo{Proof}
The preceding corollary asserts the slice-Bennequin 
inequality (with 
equality)
for $\beta$ quasipositive.  Now apply the following lemma. 
\qed 
\enddemo
\proclaim{Lemma 4 {\rm\cite{28}}} 
If the slice-Bennequin inequality holds for all 
quasipositive $\beta,$ 
then it holds for all $\beta$.
\endproclaim
\demo{Proof} Since \cite{28} is somewhat obscure, I 
resuscitate the proof.
Let
$$\beta=\s_{i_1}^{\e_1}\dotsm\s_{i_k}^{\e_k}\in B_n,\qquad
\e_j\in\{1,-1\},$$
have $p$ (resp. $\nu=k-p$) indices $j$ with $\e_j=1$
(resp. $\e_j=-1$); so $e(\beta)=p-\nu$.  If $1\le 
j_1<j_2<\dots<j_p\le k$
are the positive indices, let 
$\gamma=\s_{i_{j_1}}\negmedspace\dotsm\s_{i_{j_p}}$; so
$\gamma$ is quasipositive (in fact, positive), 
and $\chi_s(\widehat\gamma)=n-p$.  There
is a smoothly embedded surface $Q\subset S^3\times[0,1]$ 
of Euler 
characteristic $-\nu$ ($Q$ is a union of annuli with $\nu$ 
extra 
$1$-handles attached somehow) such that 
$\Bd Q\cap S^3\times\{0\}=\widehat\gamma$
and $\Bd Q\cap S^3\times\{1\}=\widehat\beta$; so  
$\vert \chi_s(\widehat\beta)-\chi_s(\widehat\gamma) \vert 
\le \nu$, 
and $\chi_s(\widehat\beta) \le n-p+\nu = n-e(\beta)$.\qed
\enddemo
\rem{Remark} Bennequin \cite{0} proved that 
$\chi(S)\le n-e(\beta)$ for all $\beta\in B_n$
and all Seifert surfaces $S$ bounded by $\widehat\beta$,
and conjectured the slice-Bennequin inequality.
\endrem

\heading 4. Nonsliceness results
\endheading
\proclaim
{Proposition} 
If the knot $K \ne O$ is strongly quasipositive, then none 
of the knots 
$D^1(K):=D(K,0,+), D^i(K):=D(D^{i-1}(K),0,+), i\ge 2,$ is 
slice.
\endproclaim
\demo{Proof}
If $K \ne O$ is strongly quasipositive, then, by Lemma~2 
and the
corollary to the Kronheimer-Mrowka Theorem, 
$D(K,0,+)$ is strongly quasipositive and not slice (because 
$\chi_s(D(K,0,+))=-1$); the proof is completed by induction.
\qed
\enddemo
\rem{Remark} Cochran and Gompf \cite{3,~Corollary~3.2} show 
that if the knot $K\ne O$ is the closure of a positive 
braid,
then $D^i(K)$ is not slice for $1\le i\le 6$.  The present 
result is 
infinitely stronger.  It would be interesting to 
understand the relation 
between being (strongly) quasipositive and ``being greater 
than 
or equal to $\Bbb T$'' in the sense of \cite{3}.
\endrem
\proclaim
{Proposition} If $p,q,r$ are all odd{\rm,} 
$\{1,-1\}\not\subset\{p,q,r\},$ and
$$
qr+rp+pq=-1,
\tag {$**$}
$$
then $\Pscr(p,q,r)$ is not slice.
\endproclaim
\rem{Remarks} (1) For $p,q,r$ odd, 
$\{1,-1\}\subset\{p,q,r\}$ iff 
$\Pscr(p,q,r)$ is an unknot, and $(**)$ iff 
$\D_{\Pscr(p,q,r)}(t)=1$.
\par
(2)~This corollary, which answers problem 1.37 in 
\cite{9}, is a special
case of results in \cite{32}.
\endrem
\demo{Proof}
Not all of $p,q,r$ have the same sign; 
we may assume $p<0<q\le r$.  By Lemma~3, if $(*)$ is true,
then $\Pscr(p,q,r)$ bounds a quasipositive Seifert surface 
of Euler
characteristic $-1$, so by \S3 it is not slice.  Suppose 
$(*)$ is false; 
then $p+q=-a, r-q=b$ with $a, b \ge 0$, so by $(**)$, 
$-1=qr+rp+pq=-(q^2+2aq+ab)$, whence
$q=1,a=0,p=-1$, and $\{1,-1\}\subset\{p,q,r\}$, contrary 
to hypothesis. 
\qed
\enddemo

\heading 5. The ``topologically locally-flat Thom 
conjecture''
\endheading

The ``Thom conjecture'' says that 
$(\vert d_a(S) \vert -1)(\vert d_a(S)\vert-2)/2 \le g(S)$ 
for any closed, oriented surface $S$ smoothly embedded in 
$\C\P^2$
of (algebraic) degree $d_a(S)$ and genus $g(S)$.  This 
conjecture is not known 
to be true, but it certainly becomes false if it
is strengthened by replacing ``smoothly embedded'' with 
``topologically 
locally-flatly embedded'' (briefly, $T$-embedded).  Let 
the {\it geometric 
degree} $d_g(S)$ of a $T$-embedded surface 
$S\subset\C\P^2$ be the minimum
number of points of intersection of a surface $S'$ 
isotopic to $S$ 
that intersects $\C\P^1_\infty$ transversally.

\rem{Claim} There is a $T$-embedded surface 
$S\subset\C\P^2$ 
with $g(S)=d_a(S)=d_g(S)=5$.
\endrem

\rem{Remark} Lee and Wilczy\'nski \cite{12} show the 
existence, 
for every $d>0$, of a $T$-embedded surface 
$W_d\subset\C\P^2$ 
with $d_a(W_d)=d$ and $g(W_d)=g_t(d)$, 
where $g_t(d)$ is the lower bound for $g(S)$
provided by classical estimates (Hsiang and Szczarba, 
Rohlin, etc.)
if $S\subset\C\P^2$ is $T$-embedded and $d_a(S)=d$;
$g_t(d)=(d-1)(d-2)/2$ for $1\le d\le 4$, and $g_t(5)=5$,
so the claim is a sharp counterexample.
The techniques of \cite{12} appear to give no control over 
$d_g(W_d)$.  It would be interesting to know if $W_d$ can 
always be taken 
to have geometric degree $d$.
\endrem
\demo{Proof {\rm(sketch)}} Follow \cite{27}; instead of
replacing a copy of $A(O\{2,3\},0)*A(O,-1)$ embedded on 
the quasipositive
Seifert surface of $O\{6,6\}$ by a $T$-embedded disk with 
the same boundary,
do the same with the copy of
$F(-3,5,7)$ embedded on the quasipositive 
Seifert surface of $O\{5,5\}$ illustrated in Fig.~2.
(By an oversight, in \cite{27} the embedding 
actually given was of $A(O\{2,3\},1)*A(O,-1)$.) \qed
\enddemo

\enddocument